\documentclass{amsart}
\usepackage{amssymb,eepic}

\newcommand\x{\mathcal X}
\newcommand\xlm{\x_{\ell,m}}
\newcommand\gf{\mathcal G}
\newcommand\CC{\mathcal C}
\newcommand\DD{\mathcal D}
\newcommand\proper{\mathfrak{P}}
\newcommand\closed{\mathfrak{C}}
\newcommand\geod{\mathfrak{G}}
\newcommand\loops{\mathfrak{L}}
\newcommand\Derives{{\buildrel * \over \Longrightarrow}}
\newcommand\tos{{\buildrel * \over \longrightarrow}}
\newcommand\free[1]{\mathbb F_{#1}}
\newcommand\set[2]{\{#1:\,#2\}}
\newcommand\C{\mathbb C}

\newcommand\Z{\mathbb Z}

\newcommand\HYP{\mathbb H}
\newcommand\cay{\operatorname{Cay}}

\newcommand\V{\operatorname{V}}

\newcommand\multifrac[2]{\begin{matrix}#1\qquad\hfill\\
  \hfill\qquad#2\end{matrix}}

\newtheorem{thm}{Theorem}[section]
\newtheorem{prop}[thm]{Proposition}
\newtheorem{lem}[thm]{Lemma} 
\newtheorem{cor}[thm]{Corollary}
\theoremstyle{definition}
\newtheorem{defn}[thm]{Definition}
\newtheorem{remark}[thm]{Remark}

\begin{document}
\title{Growth Series and Random Walks on Some Hyperbolic Graphs}
\author{Laurent Bartholdi}
\email{laurent@math.berkeley.edu}
\address{\parbox{.4\linewidth}{Department of Mathematics\\ Evans Hall\\ University of California\\ CA-94720-3840 Berkeley\\ USA}}
\author{Tullio G.\ Ceccherini--Silberstein}
\email{tceccher@mat.uniroma1.it}
\address{\parbox{.55\linewidth}{Dipartimento di Ingegneria\\ Universit\`a
  del Sannio\\ Palazzo dell'Aquila-Bosco-Lucarelli\\Corso
  Garibaldi 107, I-82100 Benevento\\ Italy}}
\date{Stockholm, December 12, 2001 --- typeset \today}
\begin{abstract}
  Consider the tessellation of the hyperbolic plane by $m$-gons, $\ell$
  per vertex. In its $1$-skeleton, we compute the growth series of
  vertices, geodesics, tuples of geodesics with common extremities. We
  also introduce and enumerate \emph{holly trees}, a family of proper
  loops in these graphs.

  We then apply Grigorchuk's result relating cogrowth and random walks
  to obtain lower estimates on the spectral radius of the Markov
  operator associated with a symmetric random walk on these graphs.
\end{abstract}
\maketitle

\section{Introduction}
We consider the graphs $\xlm$ introduced by William Floyd and Steven
Plotnick in~\cite{floyd-p:fuchsian}. These graphs are $\ell$-regular
and are the $1$-skeleton of a tessellation of the sphere (if
$(\ell-2)(m-2)<4$), the Euclidean plane (if $(\ell-2)(m-2)=4$) or the
hyperbolic plane (if $(\ell-2)(m-2)>4$) by regular $m$-gons. These
tessellations were studied by Harold Coxeter~\cite{coxeter:honey}. When
$m=\ell=4g$, then $\xlm$ is the Cayley graph of the fundamental group
$J_g = \pi_1(\Sigma_g)$ of an orientable compact surface $\Sigma_g$ of
genus $g$, with respect to the usual set of generators
$S_g=\{a_1,b_1,\dots,a_g,b_g\}$:
$$J_g = \left< a_1, b_1, \dots, a_g, b_g \left| \prod_{i=1}^g [a_i,b_i]
\right.\right>.$$

The growth series for $J_g$ with respect to $S_g$, namely
$$F_{J_g}(X) = \sum_{s\in J_g} X^{|s|} = \sum_{n=0}^\infty a_n X^n,$$
where $|s|=\min\set{t}{s=s_1\dots s_t, s_i\in S_g\cup S_g^{-1}}$
denotes the word length of $s$ with respect to $S_g$ and $a_n =
|\set{s\in J_g}{|s|=n}|$, was computed by James Cannon and Philip
Wagreich in~\cite{cannon:surfgp} and~\cite{cannon-w:surfgp} and shown
to be rational, indeed
$$F_{J_g}(X) = \frac{1 + 2X + \dots + 2X^{2g-1} + X^{2g}}{1 - (4g-2)X
  - \dots - (4g-2)X^{2g-1} + X^{2g}};$$
moreover they showed that the
denominator is a Salem polynomial. 

In~\cite{floyd-p:fuchsian} and~\cite{floyd-p:semiregular}, Floyd and
Plotnick, among other things, extended the calculations of
Cannon and Wagreich to the family $\xlm$. Fixing arbitrarily a base point
$*\in V(\xlm)$ and denoting by $|x|$ the graph distance between the
vertices $x$ and $*$, they obtained the following formula for the
growth series $F_{\ell,m}(X) = \sum_{x\in V(\xlm)} X^{|x|}$ for $m$
even:
$$F_{\ell,m}(X) = \frac{1 + 2X + \dots + 2X^{\frac m2-1} + X^{\frac m2}}{
  1 - (\ell-2)X - \dots - (\ell-2)X^{\frac m2-1} + X^{\frac m2}}$$
(they also gave a formula for $m$ odd).

In~\cite{bartholdi-c:salem} it is shown that, even in this more general
setting, the denominator is, after simplification by $(1+X)$ in case $m
\equiv 2$ mod $4$, a Salem polynomial. As a consequence of this, the
growth rates of these graphs are Salem numbers. One also obtains more
precise information about the growth coefficients: if $F_{\ell,m}(X) =
\sum_{n\ge0}a_n X^n$, then there exist constants $K$, $\lambda$ and
$R$ such that
$$K\lambda^n-R<a_n<K\lambda^n+R$$ holds for all $n$. This improves on
a result by Michel Coornaert~\cite{coornaert:mesures}.

The calculations in~\cite{cannon:surfgp}, \cite{cannon-w:surfgp}
and~\cite{floyd-p:fuchsian} are based on linear relations among
generating functions. We recover these as derivations in linear
grammars, and describe a context-free grammar enumerating more
complicated objects (see Subsection~\ref{subs:holly}).

We then compute the growth series of finite geodesics in $\xlm$
starting from $*$:
$$G_{\ell,m}(X) = \sum_{\gamma\text{ geodesic starting at $*$}}
X^{|\gamma|} = \sum_{s\in V(\xlm)} \lambda(s) X^{|s|},$$ where
$\lambda(s)$ is the number of geodesics from $*$ to $s\in V(\xlm)$,
obtaining for instance for $m$ even
$$G_{\ell,m} = \frac{1 + 2X + \dots + 2X^{\frac m2-1} + X^{\frac m2} +
  (X^{\frac m2} - X^{\frac m2-1})}{ 1 - (\ell-2)X - \dots -
  (\ell-2)X^{\frac m2-1} + X^{\frac m2} + (X^{\frac m2} - X^{\frac
    m2-1})}$$ (see the formula in Subsection~\ref{sec:ggeod} for $m$
odd); and the growth series of pairs of geodesics starting at $*$ and
with same endpoint:
$$H_{\ell,m} = \sum_{(\gamma,\delta)\,\,\parbox[t]{3.2cm}{\tiny
    geodesics starting at $*$\\ with same endpoints}}X^{|\gamma|} =
\sum_{s\in V(\xlm)} \lambda(s)^2 X^{|s|},$$
obtaining again for $m$ even
$$H_{\ell,m} = 1 + \ell\frac{X + X^2 + \dots + X^{\frac m2-1}}{ 1 -
  (\ell-2)X - \dots - (\ell-2)X^{\frac m2-1} + X^{\frac m2} + X^{\frac
    m2-1}(X-1)(2X^{\frac m2-1} - 3)}$$
(see Subsection~\ref{subs:pairs} for $m$ odd).  By computing the radius of
convergence of this last series we obtain estimates on the asymptotic
number of closed paths starting at $*$.

Better estimates in this direction are obtained by counting a larger
category of closed paths we introduced and called \emph{holly trees}; 
all these calculations are based on generatingfunctionological
methods~\cite{wilf:gf}.

Upper estimates for the spectral radius of the Markov operators
associated with a simple random walk on fundamental groups of surfaces
were calculated with various methods in~\cite{b-c-c-h:growth}. Here,
applying Grigorchuk's theorem which relates cogrowth and spectral
radius, we obtain lower estimates for the graphs $\xlm$.  We remark
that the estimates we obtain in the case of surface groups are sharper
than any other we know of (but see Remark~\ref{rem:cactus}), including
William Paschke's lower bound~\cite{paschke:norm}.  Explicitly, for
the surface group of genus $2$ we obtain the estimate of the spectral
radius $\mu_2$ of $\x_{8,8}$
$$\mu_2 \ge 0.6623.$$

The paper is organized as follows. In Section~\ref{sec:prelim} we
recall notions about graphs, groups, random walks, growth and cogrowth
as well as some facts about the $\xlm$ graphs. In
Section~\ref{sec:comput} we describe grammars and compute using them,
growth series of vertices, geodesics, pairs of geodesics and holly
trees (a family of proper closed paths) in $\xlm$. Finally in
Section~\ref{sec:lower} we use the results of Section~\ref{sec:comput}
to obtain estimates on the spectral radius of the $\xlm$.

\section{Preliminaries}\label{sec:prelim}
\subsection{Cayley graphs}
Our results are best expressed in terms of graph theory, but the
subject was first approached through group theory. We stress the
connection between these views.

\subsubsection{Groups} Let $G$ be a finitely generated group, and denote
by $S$ a finite symmetric ($S=S^{-1}$) generating subset of $G$.
Writing $S=S_+\cup S_+^{-1}$ we view $G$ as a quotient of
$\free{S_+}$, the free group on $S_+$, by a normal subgroup $N$,
namely $G = \free{S_+}/N$. In other words, if $G$ has presentation
$\langle S_+|\mathcal R\rangle$, then $N$ is the normal closure in
$\free{S_+}$ of the set of \emph{relators} $\mathcal R$.

As $G=\bigcup_{n\ge0} S^n$, we define the \emph{word length} of an
element $g\in G$ by
$$|g| = \min\set{n}{g\in S^n},$$
and the \emph{distance} between two elements $g,h\in G$ by
$$d(g,h) = |g^{-1}h|.$$
This way $(G,d)$ becomes a metric space. (If $S$ is replaced by another
finite generating set $S'$, the new distance $d'$ will be equivalent to
$d$.)

\subsubsection{Graphs} Let now $\gf$ be a graph, with vertex set $V(\gf)$
and edge set $E(\gf)\subset V(\gf)\times V(\gf)$. We suppose $E(\gf)$
is symmetric, i.e.\ is invariant under the \emph{flip} $F$ on
$V(\gf)\times V(\gf)$ defined by $F[(x,y)] = (y,x)$, and that $E(\gf)$
contains no self--loop $(x,x)$. An \emph{orientation} of $\gf$ is a
(choice of) a subset $\overset{\rightarrow}{E}\subset E(\gf)$ such
that $E(\gf)=\overset{\rightarrow}{E}\sqcup
F(\overset{\rightarrow}{E})$. If $(x,y)\in E(\gf)$ we shall refer to
the vertices $x$ and $y$ as \emph{neighbours} and write $x\sim y$. The
\emph{degree} $d(x)$ of a vertex $x$ is the number of its neighbours:
$d(x) = |\set{y}{y\sim x}|$; if this number is independent of $x$,
namely $d(x)=d(y)=d$ for all $x,y\in V(\gf)$, we say that $\gf$ is
\emph{regular of degree $d$}.

A \emph{path} in $\gf$ is a sequence $\gamma =
(\gamma_0,\dots,\gamma_n)$ of vertices of $\gf$ with $\gamma_i\sim
\gamma_{i+1}$. We call $\gamma_0$ and $\gamma_n$ respectively the
\emph{initial} and \emph{terminal} vertices of $\gamma$. The integer
$n = |\gamma|$ is the \emph{length} of $\gamma$. The \emph{inverse}
of the path $\gamma=(\gamma_0,\dots,\gamma_n)$ is
$\gamma^{-1}=(\gamma_n,\dots,\gamma_0)$; its initial vertex is the
terminal vertex of $\gamma$, and its terminal vertex is the initial
vertex of $\gamma$.

From now on we will suppose an arbitrary base point $*\in V(\gf)$ has
been fixed and will refer to $(\gf,*)$ as a rooted graph. We will also
suppose $\gf$ is \emph{connected}, i.e.\ any two points are connected
by a path, and \emph{locally finite}, i.e.\ $d(x)<\infty$ for all
$x\in V(\gf)$. (It then follows that for any $n$ there is a finite
number of paths of length $n$ starting at $*$.)

The vertex set $V(\gf)$ becomes a metric space when equipped with the
distance
$$d(x,y) = \min\set{n}{\text{there exists a path $\gamma$ with
    }\gamma_0=x, \gamma_n=y};$$ this metric is called the \emph{graph
  distance}. For $x\in V(\gf)$ we write $|x| = d(*,x)$.

A path $\gamma$ of length $n$ is said to be \emph{proper} if
$\gamma_{i-1} \neq \gamma_{i+1}$ for all $i=1,\dots,n-1$;
\emph{closed} if $\gamma_n = \gamma_0$; \emph{geodesic} if
$d(\gamma_0,\gamma_n) = n$ (and thus $d(\gamma_i,\gamma_j)=|i-j|$ for
all $i, j$); and a \emph{loop} if it is both proper and closed. We
write $\proper$, $\closed$, $\geod$, $\loops$ for the set of proper,
closed, geodesic paths and loops, respectively; and $\proper_x$,
$\closed_x$, $\geod_x$, $\loops_x$ for those starting at $x$. The
common initial and terminal vertex of a closed path is called its
\emph{root}.

The graph $\gf$ is said to be \emph{two-colourable} if there exists a
map $f:\gf\to\{\circ,\bullet\}$ such that $f(x)\neq f(y)$ whenever
$x\sim y$.

A \emph{tree} is a connected graph with no loops; in this case proper
paths are geodesics. It is obvious that all trees are two-colourable.
Any graph $\gf$ has an \emph{universal cover} $\tilde\gf$, where
$V(\tilde\gf)=\proper_*$ and $\overset{\rightarrow}{E}(\tilde\gf) =
\set{(\gamma,\delta)}{|\gamma| = |\delta|+1, \gamma_{i} = \delta_{i},
  i=0,1,\ldots |\delta|}$. It follows that $\tilde\gf$ is a tree.
There is a canonical map $\pi:\tilde\gf\to\gf$ given by $\pi(\gamma) =
\gamma_{|\gamma|}$ that is a graph homomorphism. If $\gf$ is regular,
$\tilde\gf$ is also regular, of the same degree.

If $\gf$ is a tree (or more generally a two-colourable graph) we may
orient $\gf$ by choosing for $\overset{\rightarrow}{E}$ the edges
$(x,y)$ with $x\sim y$ and $|x|<|y|$. We call this orientation the
\emph{radial orientation}.

\subsubsection{Cayley graphs}
One associates with a finitely generated group $G$ and finite
symmetric generating set $S$ (i.e., an $S \subseteq G$ with
$|S|<\infty$, $S=S^{-1}$ and $1\notin S$) its \emph{Cayley graph}
$\gf=\operatorname{Cay}(G,S)$, where $V(\gf)=G$ and $E(\gf) =
\set{(g,gs)}{g\in G, s\in S}$. The edge $e=(g,gs)$ is labelled by
$\lambda(e)=s$. The base point $*$ of $\gf$ is the vertex
corresponding to the neutral element $1$ in $G$. Cayley graphs are
connected and regular of degree $|S|$. If $G$ is free on $S_+$, its
Cayley graph is a regular tree of degree $|S|$. More generally, if $G
= \free{S_+}/N$, the universal cover $\tilde\gf$ of $\gf$ is the
Cayley graph of $\free{S_+}$ w.r.t. $S_+$, and $\pi:\tilde\gf\to\gf$
is induced by the canonical quotient map $\pi:\free{S_+}\to G$.

There is a one-to-one correspondence between paths in $\gf$ starting
at $1$ and words $w$ on the alphabet $S$; it is given explicitly by
$$\gamma=(\gamma_0,\dots,\gamma_n) \mapsto
\lambda(\gamma_0,\gamma_1)\cdots\lambda(\gamma_{n-1},\gamma_n)$$
and
$$w=w_1\cdots w_n \mapsto (1,w_1,w_1w_2,\dots,w_1\dots w_n).$$ Under
this correspondence a path $\gamma$ is proper if and only if $w_i\neq
w_{i+1}^{-1}$ for all $i$; the word corresponding to a proper path
from $1$ to~$w$ is the normal form of the element $w\in\free{S_+}$, and in
this case the correspondence is even an isometry: $|\gamma| =
|w|$. Also, $\gamma$ is geodesic if and only if $|w_1\dots
w_i| = i$ for all $i$, and $\gamma$ is closed if and only if
$\pi(w)=1$. Loops starting at $1$ correspond to elements in
$N=\ker(\pi:\free{S_+}\to G)$.

The following easy lemma is well-known:
\begin{lem}\label{lem:bicolourable}
  For a Cayley graph $\gf=\operatorname{Cay}(G,S)$ the following conditions are
  equivalent:
  \begin{enumerate}
  \item $\gf$ is two-colourable;
  \item all loops in $\gf$ have even length;
  \item all closed paths in $\gf$ have even length;
  \item all relators in $N$ have even length;
  \item there is a morphism $\phi:\,G\to \Z/2\Z=\{0,1\}$ with $\phi(S)=\{1\}$.
  \end{enumerate}
\end{lem}

\subsection{Growth and cogrowth}
Given a rooted graph $(\gf,*)$, its \emph{growth function} is
$$\gamma_\gf(n) = |\set{x\in V(\gf)}{|x|=n}|.$$
Its generating series is
$$F_\gf(X) = \sum_{n\ge0}\gamma_\gf(n)X^n.$$
The \emph{growth rate} of $\gf$ is
$$\rho_\gf = \limsup_{n\to\infty}\sqrt[n]{\gamma_\gf(n)} = R^{-1},$$
where $R$ is the radius of convergence of $F_\gf$. Note that
$\rho_\gf$ is not dependent on the choice of $*$, although $F_\gf$
is. We say $\gf$ has \emph{exponential growth} if $\rho_\gf > 1$.

If $G$ is a group generated by a finite set $S$ and $\gf$ is the
corresponding Cayley graph, the growth function of $G$ (relative to
$S$) is that of the rooted graph $(\gf,1)$, where $1$ is the unit
element in $G$.

If $Y\subseteq V(\gf)$, the \emph{relative growth} of $Y$ in $\gf$ is
$$\gamma_\gf^Y(n) = |\set{y\in Y}{|y|=n}|.$$
We define similarly its
generating series and \emph{relative growth rate}. In particular, the
relative growth rate of the set of loops $\loops_*$ viewed as a subset
of $V(\tilde\gf)$ is called the \emph{cogrowth} of $\gf$, which we
shall denote by $\alpha_\gf$; again this number $\alpha_\gf$ is not
dependent on the choice of $*$. In case $\gf=\cay(G,S)$, this is the
relative growth rate of $N=\ker(\pi)$ in $\free S$ and is called the
cogrowth of $(G,S)$.

The following estimates for a regular graph of degree $d\ge2$, due to
Rostislav Grigorchuk~\cite{grigorchuk:rw}, hold:
$$\sqrt{d-1} \le \alpha_\gf \le d-1,$$
unless $\gf$ is a tree; in that case $\alpha_\gf$ is defined to
be $1$, even though the associated growth series has infinite radius
of convergence.

%More generally, if $\phi:G\to H$ is a group epimorphism, we may
%consider the relative growth of $K=\ker(\phi)$ in $\gf=\cay(G,S)$. An
%example of this was computed by Mark Pollicott and Richard
%Sharp~\cite{pollicott:growth}, in the case $G=J_2$ and $H=(J_2)_{ab}$
%(They obtained a cogrowth series that is not rational); another
%example was given by Grigorchuk (in his thesis) for $G$ a free group and
%$H=G/[G,G]$; see also~\cite{grigorchuk-h:groups}.

\subsection{Random walks}
Let $\gf$ be a connected graph, regular of degree $d$ and denote by
$\ell^2(V(\gf)) = \{f:V(\gf) \to \C \mbox{ s.t. } \|f\|^2:=\sum_{x \in
V(\gf)} |f(x)|^2 < \infty\}$ the Hilbert space of square integrable
complex functions defined on the vertex set of $\gf$. The
\emph{simple random walk}~\cite{woess:rw} on $\gf$ is given by the stochastic
transition matrix $P=(p(x,y))_{x,y\in V(\gf)}$, with
$$p(x,y) = \begin{cases} 1/d & \text{ if } x\sim y\\ 0 & \text{
otherwise.}\end{cases}$$
The associated \emph{Markov operator} is the bounded linear operator
$M: \ell^2(V(\gf))\to\ell^2(V(\gf))$ given by
$$M(\phi)(x) = \frac1d\sum_{y\sim x}\phi(y)\qquad\text{for
  }\phi\in\ell^2(V(\gf)), x\in V(\gf).$$ Let us denote by
$\set{\delta_x}{x\in V(\gf)}$ the canonical orthonormal basis of
$\ell^2(V(\gf))$; then $p(x,y) = (M\delta_x|\delta_y)$ expresses the
probability of going to $y$ in one step, starting at $x$. More
generally $p^{(n)}(x,y) = (M^n\delta_x|\delta_y)$ is the probability
of going to $y$ in $n$ steps, starting at $x$.

The \emph{spectral radius} of the simple random walk on $\gf$ is then
$$\mu = \limsup_{n\to\infty} \sqrt[n]{p^{(n)}(*,*)} = \|M\|.$$
One can link $\mu$ to the asymptotic growth of closed paths in $\gf$:
there are approximately $(\mu d)^n$ such paths of length $n$, at least when
$n$ is even.

Harry Kesten obtained in~\cite{kesten:rwalks} the following estimates:
$$\frac{2\sqrt{d-1}}{d} \le \mu \le 1,$$
with equality on the left if
and only if $\gf$ is the regular tree of degree $d$, and equality on
the right if and only if $\gf$ is \emph{amenable}: see for
instance~\cite{dodziuk:spectra} and~\cite{ceccherini-:amen}.

\subsection{Isoperimetric constants} Let $\gf$ be a connected
graph. For a subset $K \subseteq V(\gf)$ of vertices denote by
$\partial K = \set{e = \{x,y\} \in E(\gf)}{|\{x,y\}\cap K|=1}$ the
{\it boundary} of $K$ consisting of all edges with exactly one vertex
in $K$. The number
\begin{equation}\label{iso}
  \iota(\gf) = \inf_K \frac{|\partial K|}{|K|},
\end{equation}
where the infimum is taken over all finite non-empty (connected)
subsets $K \subseteq V(\gf)$, is the \emph{(edge-)isoperimetric
  constant} of the graph $\gf$. The following fact is well-known (see
e.g.~\cite{mohar:isoperimetric,mohar-w:spectra}):
\begin{thm} Let $\gf$ be a regular graph of degree $d$. Denoting by
  $\iota = \iota(\gf)$ and $\mu$ the isoperimetric constant and the
  spectral radius of the simple random walk on $\gf$, respectively,
  one has
%\begin{equation}\label{iotamu}
%\frac{d^2}{d-1}(1-\mu) \leq \iota \leq d \sqrt{1 - \mu^2}
%\end{equation} and, equivalently,
  \begin{equation}\label{muiota}
   1 - \left(\frac{d-1}{d^2}\right)\iota \leq \mu \leq \sqrt{1 - \left(\frac{\iota}{d}\right)^2}.
  \end{equation}
\end{thm}

Note that, by Kesten's results and (\ref{muiota}), $\gf$ is amenable
if and only if $\iota(\gf) = 0$ (see \cite{ceccherini-:amen} for more
on this).

\subsection{The Grigorchuk formula}
In his thesis, Grigorchuk~\cite{grigorchuk:rw} found the following relation
between the cogrowth $\alpha$ of $\gf$ and the spectral radius $\mu$
of the simple random walk on a regular graph $\gf$ of degree $d$,
namely
\begin{equation}\label{eq:grigorchuk}
  \mu = \begin{cases}
   \frac{2\sqrt{d-1}}d & \text{ if } \alpha\le\sqrt{d-1} \\
   \frac{\sqrt{d-1}}d\left(\frac{\sqrt{d-1}}{\alpha} +
    \frac{\alpha}{\sqrt{d-1}}\right) & \text{ if } \alpha>\sqrt{2d-1}.
  \end{cases}
\end{equation}
A generalization of~(\ref{eq:grigorchuk}), along with a simple proof
can be found in~\cite{bartholdi:cogrowth}.

\begin{remark}\label{rem:alphamu}
  Since the function $\mu(t) = (\sqrt{d-1}/d)\left(\sqrt{d-1}/t +
    t/\sqrt{d-1}\right)$ is an increasing function of $t$ for
  $t\in[\sqrt{d-1},d-1]$, any lower bound on $\alpha$ yields a lower
  bound on $\mu$.
\end{remark}

\subsection{Surface groups}
Let $J_g = \pi_1(\Sigma_g)$ denote the fundamental group of an
orientable surface of genus $g$. This group is finitely generated and
admits the following presentation:
$$J_g = \left< a_1, b_1, \dots, a_g, b_g \left| \prod_{i=1}^g [a_i,b_i]
\right.\right>.$$

Similarly denoting by $J'_g = \pi(\Sigma'_g)$ the fundamental group of
a non-orientable surface of genus $g\ge1$, we have
$$J'_g = \left< a_1, \dots, a_g \left| \prod_{i=1}^g a_i^2
\right.\right>.$$
A simple justification of these presentations is given in the first
chapters
of~\cite{massey:at}.

The Cayley graphs corresponding to these groups with generators
specified as above constitute the 1-skeleton of planar tessellations,
the underlying space being the sphere $\mathbb S^2$ for $J_0$ and
$J'_1$; the Euclidean plane $\mathbb R^2$ for $J_1$ and $J'_2$; and
the hyperbolic plane $\mathbb H^2$ for $J_g$, $g\ge2$ and $J'_g$,
$g\ge3$. The
2-cells of the tessellation are $4g$-gons (respectively $2g$-gons) for
$J_g$ (respectively $J'_g$).
These facts are used and developed in~\cite{cannon:surfgp}
and~\cite[chapter~III, \S 5]{lyndon-s:cgt}.

\subsection{The graphs $\xlm$}\label{subs:graphsxlm}
We shall be concerned with the more general graphs $\xlm$ that are the
1-skeleton of regular tessellations of a constant-curvature surface
consisting of $m$-gons, with vertex angles
$2\pi/\ell$~\cite{floyd-p:fuchsian}. We assume that $\ell\ge3$ and
$m\ge3$. It is proved for instance in~\cite{coxeter:honey} that such a
tiling by polygons exists and is unique up to isometry. We fix once
and for all an arbitrary base point $*\in V(\xlm)$.

When $(\ell-2)(m-2) > 4$, these graphs can be embedded
quasi-isometrically in $\HYP^2$, the Poincar\'e-Lobachevsky plane.
When $(\ell-2)(m-2) = 4$, they embed in Euclidean plane $\mathbb R^2$.
When $(\ell-2)(m-2) < 4$, they embed in $\mathbb S^2$, the sphere.
For instance, if $(\ell,m)=(3,3)$ we have the tetrahedral tessellation
of the sphere; if $(\ell,m)=(4,3)$ the octahedral tessellation; if
$(\ell,m)=(5,3)$ the icosahedral tessellation; if $(\ell,m)=(3,4)$ the
cubical tessellation; and if $(\ell,m)=(3,5)$ the dodecahedral
tessellation of the sphere. For $(\ell,m)=(3,6), (4,4), (6,3)$ we have
respectively the hexagonal, square and triangular Euclidean
tessellations; and in all other cases a hyperbolic tessellation.
When $\ell=m=4g$ they are the
Cayley graphs of $J_g$ and when $\ell=m=2g$ the Cayley graphs of
$J'_g$, with respect to their canonical generating sets.

By Lemma~\ref{lem:bicolourable}, these graphs are two-colourable if
and only if $m$ is even.

Thanks to the underlying space and the presentation of these graphs as
tessellations, there is a natural notion of \emph{cell} or \emph{face}.
These are closures of connected components of the complement of $\xlm$
in the embedding space. Cell boundaries are precisely the simple
loops in $\xlm$ of length $m$.

Recall the graph-theoretical notion of \emph{dual graph}. Given a
graph $\gf=(V,E)$ embedded in a surface with set of cells $\mathcal
C$, its dual is the graph $\check \gf=(\mathcal C,F)$, with an edge
$(C,D)\in F$ whenever $C$ and $D$, viewed as two-cells in $\gf$, share
a common edge in $E$. In our setting the dual of $\xlm$ is
$\x_{m,\ell}$. In particular, the Cayley graphs of surface groups are
self-dual.

\begin{defn}\label{def:cone}
  Let $\gf$ be a rooted graph. The \emph{cone} of a vertex $x\in
  V(\gf)$ is the subgraph spanned by the set of vertices
  $$\set{y\in V(\gf)}{\exists\gamma\in \geod_*\text{ with }
   \gamma_i=x, \gamma_j=y, \text{ for some }i\le j}.$$ This is the
  set of vertices which may be joined starting from the base point $*$ by a
  geodesic passing through $x$, namely $\gamma = (*, \ldots, x, \ldots, y)$.
  
  Two cones are said to be \emph{equivalent} if there exists a
  rooted-graph-isomorphism mapping one cone onto the other.  The
  \emph{cone type} of a vertex $x$ is the class $t(x)$ of its cone
  modulo this equivalence.
\end{defn}

Call a vertex $x$ a \emph{successor} of $y$ if $x\sim y$ and $|x|>|y|$
(hence $|x| = |y|+1$). Symmetrically $y$ is a \emph{predecessor} of
$x$. If $x\sim y$ and $|x|=|y|$ call $x$ and $y$ \emph{peers}.
Graphically, if $y$ is a successor of $x$, we orient the edge $(x,y)$
in the direction $x\to y$ and we say that this arrow \emph{exits} from
$x$ and \emph{enters} in $y$. It follows easily from the definition
that the cone type $t$ of $x$ determines the number $\ell(t)$ and the
cone types of its successors.

We label the edges of $\xlm$ using the alphabet
$\{e_1,\dots,e_\ell,e_p\}$, as follows: arbitrarily label clockwise
$e_1,\dots,e_\ell$ the edges exiting from $*$. At each vertex
$x\in\V(\xlm)$, label from left to right $e_1,\dots,e_{\ell(t)}$ the
edges exiting from $x$. If $x$ has a peer $y$, label $e_p$ the edge
between $x$ and $y$.

Note that if $x$ and $y$ have the same cone type, then there is an
edge-label-preserving isomorphism between their cones.

The graphs $\xlm$ have a finite number of cone types. This was first
proven by Cannon for the Cayley graphs of surface
groups~\cite{cannon:surfgp} and later extended to the hyperbolic
$\xlm$ by Floyd and Plotnick~\cite{floyd-p:fuchsian}. The cone types
can be described in several equivalent ways among which we propose the
following. First suppose $m$ is even.
\begin{defn}\label{defn:type}
  For $x\in\V(\xlm)$, with $m$ even, set
  $$t(x) = \max_{(F,y):\,F\ni x, \ y\in F}\left(|x|-|y|\right) = |x| - \min_{(F,y):\,F\ni x,\ y\in F}|y|.$$
In the expression above $F$ denotes a cell of the tessellation, i.e.\ an
  $m$-gon.
\end{defn}
One has $t(x)\in\{0,\dots,m/2\}$. Also, $t(x)=0$ holds only for $x=*$.

Furthermore, cone types appear in two mirror versions, depending on
which side (left or right) contains the $2$-cell closest to origin.
We extend the cone types of vertices in $\xlm$ as follows: if $v$ has
type $0$ or $1$, then its extended type is the same. If $v$ has type
$t\ge2$, and $v$ is on the left side of its lowest adjacent $2$-cell,
then its \emph{extended type} is $(L,t)$. If $v$ is on the right side
of that $2$-cell, its extended type is $(R,t)$.

Extended types will be used in the computations of
Section~\ref{sec:comput}. For now, we describe diagrammatically the
(left) cone types, letting $t$ be any type in $\{1,\dots,m/2\}$:
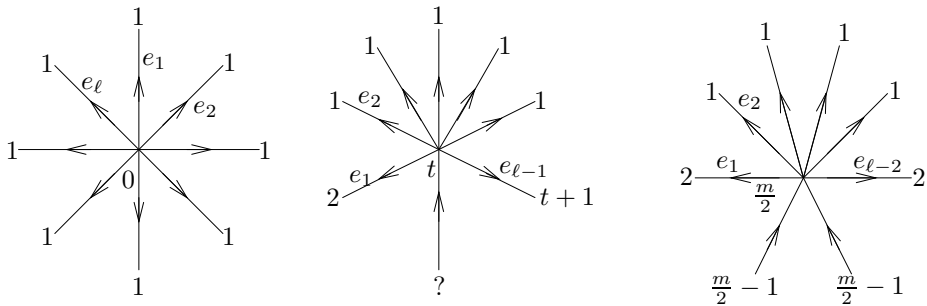
\[\setlength{\unitlength}{0.00083333in}
\begingroup\makeatletter\ifx\SetFigFont\undefined%
\gdef\SetFigFont#1#2#3#4#5{%
  \reset@font\fontsize{#1}{#2pt}%
  \fontfamily{#3}\fontseries{#4}\fontshape{#5}%
  \selectfont}%
\fi\endgroup%
{
\begin{picture}(1670,1836)(0,-10)
\put(850,1434){\makebox(0,0)[lb]{\smash{$e_1$}}}
\put(1170,1134){\makebox(0,0)[lb]{\smash{$e_2$}}}
\put(470,1284){\makebox(0,0)[lb]{\smash{$e_\ell$}}}
\path(831,909)(1281,909)
\path(1161.000,879.000)(1281.000,909.000)(1161.000,939.000)
\path(831,909)(1131,1209)
\path(1067.360,1102.934)(1131.000,1209.000)(1024.934,1145.360)
\path(831,909)(531,1209)
\path(637.066,1145.360)(531.000,1209.000)(594.640,1102.934)
\path(831,909)(831,1359)
\path(861.000,1239.000)(831.000,1359.000)(801.000,1239.000)
\path(594.640,715.066)(531.000,609.000)(637.066,672.640)
\path(531,609)(831,909)
\path(1024.934,672.640)(1131.000,609.000)(1067.360,715.066)
\path(1131,609)(831,909)
\path(831,909)(381,909)
\path(501.000,939.000)(381.000,909.000)(501.000,879.000)
\path(831,909)(831,459)
\path(801.000,579.000)(831.000,459.000)(861.000,579.000)
\path(81,909)(531,909)
\path(306,384)(681,759)
\path(306,1434)(681,1059)
\path(831,1659)(831,1209)
\path(1356,1434)(981,1059)
\path(1581,909)(1131,909)
\path(1356,384)(981,759)
\path(831,159)(831,609)
\put(786,1686){\makebox(0,0)[lb]{\smash{$1$}}}
\put(1350,312){\makebox(0,0)[lb]{\smash{$1$}}}
\put(786,0){\makebox(0,0)[lb]{\smash{$1$}}}
\put(0,849){\makebox(0,0)[lb]{\smash{$1$}}}
\put(216,1395){\makebox(0,0)[lb]{\smash{$1$}}}
\put(1356,1401){\makebox(0,0)[lb]{\smash{$1$}}}
\put(219,312){\makebox(0,0)[lb]{\smash{$1$}}}
\put(1578,846){\makebox(0,0)[lb]{\smash{$1$}}}
\put(726,669){\makebox(0,0)[lb]{\smash{$0$}}}
\end{picture}
}
\qquad
\setlength{\unitlength}{0.00083333in}
\begingroup\makeatletter\ifx\SetFigFont\undefined%
\gdef\SetFigFont#1#2#3#4#5{%
  \reset@font\fontsize{#1}{#2pt}%
  \fontfamily{#3}\fontseries{#4}\fontshape{#5}%
  \selectfont}%
\fi\endgroup%
{
\begin{picture}(1876,1842)(0,-10)
\put(142,704){\makebox(0,0)[lb]{\smash{$e_1$}}}
\put(197,1199){\makebox(0,0)[lb]{\smash{$e_2$}}}
\put(1077,759){\makebox(0,0)[lb]{\smash{$e_{\ell-1}$}}}
\path(702,909)(927,1284)
\path(890.985,1165.666)(927.000,1284.000)(839.536,1196.536)
\path(702,909)(330,726)
\path(424.434,805.889)(330.000,726.000)(450.919,752.051)
\path(702,909)(1077,723)
\path(956.167,749.446)(1077.000,723.000)(982.828,803.197)
\path(702,909)(330,1095)
\path(450.748,1068.167)(330.000,1095.000)(423.915,1014.502)
\path(702,909)(1092,1104)
\path(998.085,1023.502)(1092.000,1104.000)(971.252,1077.167)
\path(702,909)(702,1359)
\path(732.000,1239.000)(702.000,1359.000)(672.000,1239.000)
\path(702,159)(702,642)
\path(732.000,522.000)(702.000,642.000)(672.000,522.000)
\path(702,459)(702,909)
\path(999,759)(1302,609)
\path(1002,1059)(1302,1209)
\path(102,609)(402,759)
\path(102,1209)(402,1059)
\path(702,1209)(702,1659)
\path(702,909)(477,1284)
\path(564.464,1196.536)(477.000,1284.000)(513.015,1165.666)
\path(591,1092)(318,1545)
\path(816,1092)(1074,1539)
\put(627,738){\makebox(0,0)[lb]{\smash{$t$}}}
\put(666,0){\makebox(0,0)[lb]{\smash{$?$}}}
\put(657,1692){\makebox(0,0)[lb]{\smash{$1$}}}
\put(1074,1506){\makebox(0,0)[lb]{\smash{$1$}}}
\put(231,1506){\makebox(0,0)[lb]{\smash{$1$}}}
\put(12,1158){\makebox(0,0)[lb]{\smash{$1$}}}
\put(1296,1158){\makebox(0,0)[lb]{\smash{$1$}}}
\put(1332,561){\makebox(0,0)[lb]{\smash{$t+1$}}}
\put(0,558){\makebox(0,0)[lb]{\smash{$2$}}}
\end{picture}
}
\qquad
\setlength{\unitlength}{0.00083333in}
\begingroup\makeatletter\ifx\SetFigFont\undefined%
\gdef\SetFigFont#1#2#3#4#5{%
  \reset@font\fontsize{#1}{#2pt}%
  \fontfamily{#3}\fontseries{#4}\fontshape{#5}%
  \selectfont}%
\fi\endgroup%
{
\begin{picture}(1686,1740)(0,-10)
\put(226,792){\makebox(0,0)[lb]{\smash{$e_1$}}}
\put(361,1182){\makebox(0,0)[lb]{\smash{$e_2$}}}
\put(1081,792){\makebox(0,0)[lb]{\smash{$e_{\ell-2}$}}}
\path(771,732)(1221,732)
\path(1101.000,702.000)(1221.000,732.000)(1101.000,762.000)
\path(771,732)(1146,1107)
\path(1082.360,1000.934)(1146.000,1107.000)(1039.934,1043.360)
\path(771,732)(396,1107)
\path(502.066,1043.360)(396.000,1107.000)(459.640,1000.934)
\path(771,732)(321,732)
\path(441.000,762.000)(321.000,732.000)(441.000,702.000)
\path(771,732)(912,1257)
\path(909.848,1133.326)(912.000,1257.000)(851.901,1148.888)
\path(96,732)(621,732)
\path(1446,732)(921,732)
\path(771,732)(627,1260)
\path(687.517,1152.122)(627.000,1260.000)(629.631,1136.335)
\path(471,132)(621,432)
\path(594.167,311.252)(621.000,432.000)(540.502,338.085)
\path(1071,132)(921,432)
\path(1001.498,338.085)(921.000,432.000)(947.833,311.252)
\path(546,282)(771,732)
\path(996,282)(771,732)
\path(771,732)(996,1557)
\path(246,1257)(621,882)
\path(546,1557)(771,732)
\path(921,882)(1296,1257)
\put(455,576){\makebox(0,0)[lb]{\smash{$\frac m2$}}}
\put(180,0){\makebox(0,0)[lb]{\smash{$\frac m2-1$}}}
\put(156,1215){\makebox(0,0)[lb]{\smash{$1$}}}
\put(495,1590){\makebox(0,0)[lb]{\smash{$1$}}}
\put(984,1587){\makebox(0,0)[lb]{\smash{$1$}}}
\put(1302,1221){\makebox(0,0)[lb]{\smash{$1$}}}
\put(963,3){\makebox(0,0)[lb]{\smash{$\frac m2-1$}}}
\put(0,687){\makebox(0,0)[lb]{\smash{$2$}}}
\put(1449,681){\makebox(0,0)[lb]{\smash{$2$}}}
\end{picture}
}
\]

\begin{table}
\begin{center}
\begin{tabular}{r|c@{\quad}c@{ }c@{ }c}
vertex & \multicolumn{3}{c}{successors of type}& number of \\
type $t$     & $1$       & $2$  & $t+1$ & predecessors
\\[2pt]\hline
\rule{0pt}{2.5ex}
$0$         & $\ell$     & $0$  & $0$  & $0$ \\[1pt]
$0<t<\frac m2$ & $\ell-3$    & $1$  & $1$  & $1$ \\[1pt]
$\frac m2$    & $\ell-4$    & $2$  & $0$  & $2$ \\[1ex]
\end{tabular}
\end{center}
\caption{The cone type data for even $m$}\label{table:cte}
\end{table}

The data we will use are summarized in Table~\ref{table:cte}. For
instance, a vertex of type $m/2$ has two successors of type $2$,
$\ell-4$ of type $1$ and two predecessors of type $m/2-1$.

We now consider the case when $m$ is odd. As said earlier $\xlm$ is
not two-colourable: in fact there exist peers in $\xlm$. If we
suppress all edges $(x,y)$ connecting peers we obtain a tessellation by
$(2m-2)$-gons, which we call temporarily $\widehat\xlm$ (note that the
corresponding graph is no more regular).
\begin{defn}
  For $x\in\V(\xlm)$, with $m$ odd, set
  $$t(x) = \max_{(F,y):\,F\ni x, y\in F}\left(|x|-|y|\right) = |x| - \min_{(F,y):\,F\ni x,\ y\in F}|y|.$$
In the expression above $F$ denotes a cell of $\widehat\xlm$, i.e.\ a
$(2m-2)$-gon.
\end{defn}
Obviously $t(x)\in\{0,\dots,m-1\}$. Also, $t(x)=0$ holds only for $x=*$.

As before, the (left) cone types can be represented by
\begin{align*}
\setlength{\unitlength}{0.00083333in}
\begingroup\makeatletter\ifx\SetFigFont\undefined%
\gdef\SetFigFont#1#2#3#4#5{%
  \reset@font\fontsize{#1}{#2pt}%
  \fontfamily{#3}\fontseries{#4}\fontshape{#5}%
  \selectfont}%
\fi\endgroup%
{\begin{picture}(1670,1836)(0,-10)
\put(850,1434){\makebox(0,0)[lb]{\smash{$e_1$}}}
\put(1170,1134){\makebox(0,0)[lb]{\smash{$e_2$}}}
\put(470,1284){\makebox(0,0)[lb]{\smash{$e_\ell$}}}
\path(831,909)(1281,909)
\path(1161.000,879.000)(1281.000,909.000)(1161.000,939.000)
\path(831,909)(1131,1209)
\path(1067.360,1102.934)(1131.000,1209.000)(1024.934,1145.360)
\path(831,909)(531,1209)
\path(637.066,1145.360)(531.000,1209.000)(594.640,1102.934)
\path(831,909)(831,1359)
\path(861.000,1239.000)(831.000,1359.000)(801.000,1239.000)
\path(594.640,715.066)(531.000,609.000)(637.066,672.640)
\path(531,609)(831,909)
\path(1024.934,672.640)(1131.000,609.000)(1067.360,715.066)
\path(1131,609)(831,909)
\path(831,909)(381,909)
\path(501.000,939.000)(381.000,909.000)(501.000,879.000)
\path(831,909)(831,459)
\path(801.000,579.000)(831.000,459.000)(861.000,579.000)
\path(81,909)(531,909)
\path(306,384)(681,759)
\path(306,1434)(681,1059)
\path(831,1659)(831,1209)
\path(1356,1434)(981,1059)
\path(1581,909)(1131,909)
\path(1356,384)(981,759)
\path(831,159)(831,609)
\put(786,1686){\makebox(0,0)[lb]{\smash{$1$}}}
\put(1350,312){\makebox(0,0)[lb]{\smash{$1$}}}
\put(786,0){\makebox(0,0)[lb]{\smash{$1$}}}
\put(0,849){\makebox(0,0)[lb]{\smash{$1$}}}
\put(216,1395){\makebox(0,0)[lb]{\smash{$1$}}}
\put(1356,1401){\makebox(0,0)[lb]{\smash{$1$}}}
\put(219,312){\makebox(0,0)[lb]{\smash{$1$}}}
\put(1578,846){\makebox(0,0)[lb]{\smash{$1$}}}
\put(726,669){\makebox(0,0)[lb]{\smash{$0$}}}
\end{picture}
}
&\qquad
\setlength{\unitlength}{0.00083333in}
\begingroup\makeatletter\ifx\SetFigFont\undefined%
\gdef\SetFigFont#1#2#3#4#5{%
  \reset@font\fontsize{#1}{#2pt}%
  \fontfamily{#3}\fontseries{#4}\fontshape{#5}%
  \selectfont}%
\fi\endgroup%
{\begin{picture}(1876,1842)(0,-10)
\put(122,704){\makebox(0,0)[lb]{\smash{$e_1$}}}
\put(177,1209){\makebox(0,0)[lb]{\smash{$e_2$}}}
\put(1077,759){\makebox(0,0)[lb]{\smash{$e_{\ell-1}$}}}
\path(702,909)(927,1284)
\path(890.985,1165.666)(927.000,1284.000)(839.536,1196.536)
\path(702,909)(330,726)
\path(424.434,805.889)(330.000,726.000)(450.919,752.051)
\path(702,909)(1077,723)
\path(956.167,749.446)(1077.000,723.000)(982.828,803.197)
\path(702,909)(330,1095)
\path(450.748,1068.167)(330.000,1095.000)(423.915,1014.502)
\path(702,909)(1092,1104)
\path(998.085,1023.502)(1092.000,1104.000)(971.252,1077.167)
\path(702,909)(702,1359)
\path(732.000,1239.000)(702.000,1359.000)(672.000,1239.000)
\path(702,159)(702,642)
\path(732.000,522.000)(702.000,642.000)(672.000,522.000)
\path(702,459)(702,909)
\path(999,759)(1302,609)
\path(1002,1059)(1302,1209)
\path(102,609)(402,759)
\path(102,1209)(402,1059)
\path(702,1209)(702,1659)
\path(702,909)(477,1284)
\path(564.464,1196.536)(477.000,1284.000)(513.015,1165.666)
\path(591,1092)(318,1545)
\path(816,1092)(1074,1539)
\put(627,738){\makebox(0,0)[lb]{\smash{$t$}}}
\put(666,0){\makebox(0,0)[lb]{\smash{$?$}}}
\put(657,1692){\makebox(0,0)[lb]{\smash{$1$}}}
\put(1074,1506){\makebox(0,0)[lb]{\smash{$1$}}}
\put(231,1506){\makebox(0,0)[lb]{\smash{$1$}}}
\put(12,1158){\makebox(0,0)[lb]{\smash{$1$}}}
\put(1296,1158){\makebox(0,0)[lb]{\smash{$1$}}}
\put(1332,561){\makebox(0,0)[lb]{\smash{$t+1$}}}
\put(0,558){\makebox(0,0)[lb]{\smash{$2$}}}
\end{picture}
}
\\
\setlength{\unitlength}{0.00083333in}
\begingroup\makeatletter\ifx\SetFigFont\undefined%
\gdef\SetFigFont#1#2#3#4#5{%
  \reset@font\fontsize{#1}{#2pt}%
  \fontfamily{#3}\fontseries{#4}\fontshape{#5}%
  \selectfont}%
\fi\endgroup%
{\begin{picture}(2160,1914)(0,-10)
\put(272,721){\makebox(0,0)[lb]{\smash{$e_1$}}}
\put(242,1121){\makebox(0,0)[lb]{\smash{$e_2$}}}
\put(1227,956){\makebox(0,0)[lb]{\smash{$e_p$}}}
\path(822,906)(1032,1281)
\path(999.543,1161.641)(1032.000,1281.000)(947.192,1190.957)
\path(822,906)(594,1245)
\path(685.864,1162.168)(594.000,1245.000)(636.077,1128.683)
\path(822,906)(447,720)
\path(541.172,800.197)(447.000,720.000)(567.833,746.446)
\path(822,906)(426,1038)
\path(549.329,1028.513)(426.000,1038.000)(530.355,971.592)
\path(822,906)(1206,1098)
\path(1112.085,1017.502)(1206.000,1098.000)(1085.252,1071.167)
\path(822,906)(1497,906)
\path(822,156)(822,672)
\path(852.000,552.000)(822.000,672.000)(792.000,552.000)
\path(255,621)(522,756)
\path(1422,1206)(1122,1056)
\path(822,306)(822,906)
\path(822,906)(423,1494)
\path(1122,981)(1122,831)
\path(1161,981)(1161,831)
\path(822,906)(1176,1545)
\path(822,906)(147,1131)
\path(822,906)(825,1734)
\path(822,906)(822,1356)
\path(852.000,1236.000)(822.000,1356.000)(792.000,1236.000)
\put(786,0){\makebox(0,0)[lb]{\smash{$?$}}}
\put(60,1080){\makebox(0,0)[lb]{\smash{$1$}}}
\put(1422,1149){\makebox(0,0)[lb]{\smash{$1$}}}
\put(783,1764){\makebox(0,0)[lb]{\smash{$1$}}}
\put(1170,1512){\makebox(0,0)[lb]{\smash{$1$}}}
\put(336,1512){\makebox(0,0)[lb]{\smash{$1$}}}
\put(1530,861){\makebox(0,0)[lb]{\smash{$\frac{m-1}2$}}}
\put(852,696){\makebox(0,0)[lb]{\smash{$\frac{m-1}2$}}}
\put(0,501){\makebox(0,0)[lb]{\smash{$\frac{m+1}2$}}}
\end{picture}
}
&\qquad
\setlength{\unitlength}{0.00083333in}
\begingroup\makeatletter\ifx\SetFigFont\undefined%
\gdef\SetFigFont#1#2#3#4#5{%
  \reset@font\fontsize{#1}{#2pt}%
  \fontfamily{#3}\fontseries{#4}\fontshape{#5}%
  \selectfont}%
\fi\endgroup%
{\begin{picture}(1830,1725)(0,-10)
\put(177,767){\makebox(0,0)[lb]{\smash{$e_1$}}}
\put(347,1147){\makebox(0,0)[lb]{\smash{$e_2$}}}
\put(1077,767){\makebox(0,0)[lb]{\smash{$e_{\ell-2}$}}}
\path(777,717)(1227,717)
\path(1107.000,687.000)(1227.000,717.000)(1107.000,747.000)
\path(777,717)(1152,1092)
\path(1088.360,985.934)(1152.000,1092.000)(1045.934,1028.360)
\path(777,717)(402,1092)
\path(508.066,1028.360)(402.000,1092.000)(465.640,985.934)
\path(777,717)(327,717)
\path(447.000,747.000)(327.000,717.000)(447.000,687.000)
\path(777,717)(918,1242)
\path(915.848,1118.326)(918.000,1242.000)(857.901,1133.888)
\path(102,717)(627,717)
\path(1452,717)(927,717)
\path(777,717)(633,1245)
\path(693.517,1137.122)(633.000,1245.000)(635.631,1121.335)
\path(477,117)(627,417)
\path(600.167,296.252)(627.000,417.000)(546.502,323.085)
\path(1077,117)(927,417)
\path(1007.498,323.085)(927.000,417.000)(953.833,296.252)
\path(552,267)(777,717)
\path(1002,267)(777,717)
\path(777,717)(1002,1542)
\path(252,1242)(627,867)
\path(552,1542)(777,717)
\path(927,867)(1302,1242)
\put(162,1200){\makebox(0,0)[lb]{\smash{$1$}}}
\put(501,1575){\makebox(0,0)[lb]{\smash{$1$}}}
\put(990,1572){\makebox(0,0)[lb]{\smash{$1$}}}
\put(1308,1206){\makebox(0,0)[lb]{\smash{$1$}}}
\put(936,6){\makebox(0,0)[lb]{\smash{$m-2$}}}
\put(315,0){\makebox(0,0)[lb]{\smash{$m-2$}}}
\put(520,555){\makebox(0,0)[lb]{\smash{$m-1$}}}
\put(1458,666){\makebox(0,0)[lb]{\smash{$2$}}}
\put(0,669){\makebox(0,0)[lb]{\smash{$2$}}}
\end{picture}
}
\end{align*}

\begin{table}
\begin{center}
\begin{tabular}{r|c@{\quad}c@{ }c@{ }c@{ }c}
vertex & \multicolumn{3}{c}{successors of type}& number of   & number
of \\
type $t$     & $1$  & $2$  & $t+1$ & predecessors & peers \\[2pt]
\hline
\rule{0pt}{2.5ex}
$0$         & $\ell$     & $0$  & $0$  & $0$  & $0$ \\[1pt]
$0<t<m-1, t\neq\frac{m-1}{2}$
           & $\ell-3$    & $1$  & $1$  & $1$  & $0$ \\[1pt]
$\frac{m-1}{2}$ & $\ell-4$    & $1$  & $1$  & $1$  & $1$ \\[1pt]
$m-1$       & $\ell-4$    & $2$  & $0$  & $2$  & $0$ \\[1ex]
\end{tabular}
\end{center}
\caption{The cone type data for odd $m$}\label{table:cto}
\end{table}

The data we will use are summarized in Table~\ref{table:cto}.
The two predecessors of a vertex of type $m-1$ are of type $m-2$; the peer
of a vertex of type $(m-1)/2$ also has type $(m-1)/2$.

\subsection{Hyperbolicity of the graphs $\xlm$}
Let $(X,d)$ be a geodesic metric space and $\delta >0$. We say $X$ is
Gromov-$\delta$-\emph{hyperbolic} if \emph{geodesic triangles} are
{$\delta$-thin} in the following sense: for all $x, y, z \in X$ one
has
\begin{equation}\label{eq:triangle}
\forall t \in \overline{yz}:\; d(t, \overline{xy} \cup \overline{xz}) < \delta
\end{equation}
where $\overline{yz}$ is a geodesic segment between $y$ and $z$.

As a special case, two geodesics $\gamma$ and $\gamma'$ from
$x$ and $y$ in $X$ are at distance at most $\delta$ apart:
\begin{equation}\label{eq:hyp}
  \forall t \in \gamma,\text{ there exists }t' \in \gamma'\text{ such
   that }d(t,t') < \delta
\end{equation}
(take $z = y$ in~(\ref{eq:triangle})).

A finitely generated group $G$ is \emph{hyperbolic} if its Cayley
graph $\gf = \operatorname{Cay}(G,S)$ is hyperbolic as a graph (with a
constant of hyperbolicity $\delta_S$ which depends on the generating
set $S$): this is well defined since for any other generating
set $S'$ the graph $\gf' = \operatorname{Cay}(G,S')$ will still be
hyperbolic (with a possibly different constant of hyperbolicity
$\delta_{S'}$). 

Finite graphs and finite groups are clearly hyperbolic.

The groups $J_g$ are hyperbolic~\cite{ghys-h:gromov,gromov:asympt} as
soon as $g\neq 1$, and $J'_g$ are hyperbolic as soon as $g\neq 1, 2$. The
other cases are well-known: $J_0$ is trivial, $J_1$ is $\Z^2$, $J'_2$
is the two-group and $J'_3$ is a twofold extension of $\Z^2$.
 
Our graphs $\xlm$ are hyperbolic whenever $(\ell-2)(m-2) \neq 4$:
indeed one may take $\delta = w$ for even $m=2w$ and odd $m=2w+1$.

The following lemma strengthens for the graphs $\xlm$ the situation
described in~(\ref{eq:hyp}):
\begin{lem}\label{lem:hyper}
  Let $\gf$ be a hyperbolic $\xlm$; let $\gamma$ and $\gamma'$ be two
  geodesics in $\gf$ of length $n$, with same extremities $\gamma_0 =
  \gamma'_0$ and $\gamma_n = \gamma'_n$. Then for all
  $i\in\{0,\dots,n\}$, there exists a cell $F_i$ in $\gf$ such that
  both $\gamma_i$ and $\gamma'_i$ belong to $F_i$.
\end{lem}
\begin{proof}
  The proof is not difficult and relies on Euler characteristic
  considerations; see~\cite{b-c-c-h:growth} and~\cite{zuk:norm}.
\end{proof}

\section{Recursive Computations}\label{sec:comput}
We present in this section the methods used to compute the growth
series of various objects related to $\xlm$: vertices, geodesics, etc.
The unifying notion is that of \emph{production grammar}. The
procedure we will follow for each family of objects is: construct a
grammar and compute the growth series of its associated language.

\subsection{Grammars}
A \emph{context-free grammar}~\cite{harrison:flt} is a tuple
$\Gamma=(N,T,S,R)$, with $N$ and $T$ disjoint finite sets, called
respectively the \emph{non-terminal} and \emph{terminal} alphabets;
$S\in N$ called the \emph{axiom}; and $R$ a finite subset of $N\times
(N\cup T)^*$ called the \emph{set of (production) rules}, where for an alphabet
$A$, we denote by $A^*$ the set of finite-length words over $A$. A
rule $(X,w)\in R$ is conveniently written $X\to w$.

Let $\Gamma=(N,T,S,R)$ be a context-free grammar, and $u,v\in(N\cup
T)^*$ two words. We say $v$ \emph{is derived from} $u$, and write $u
\Longrightarrow v$, if there is a rule $(X,w)\in R$ and factorizations
$u=aXb$, $v=awb$ for some $a\in(N\cup T)^*$ and $b\in T^*$. Let
$\Derives$ be the transitive closure of $\Longrightarrow$: one has
$u\Derives v$ precisely when there is a sequence $u=u_0
\Longrightarrow u_1 \Longrightarrow\dots \Longrightarrow u_n=v$. The
\emph{language} of $\Gamma$ is
$$L(\Gamma)=\set{w\in T^*}{S\Derives w}.$$
The \emph{intermediate derivations} of $\Gamma$ are the $w\in(N\cup
T)^*$ such that $S\Derives w$.

If all rules $(X,w)\in R$ satisfy $w\in T^*N$, the grammar is called
\emph{right-linear}; if they all satisfy $w\in NT^*$, the grammar is
called \emph{left-linear}. In both cases the grammar and associated
language are also termed \emph{regular}. If for any word $w\in
L(\Gamma)$ there is a unique sequence of derivations that yield $w$
from $S$, the grammar is called \emph{unambiguous}. The \emph{growth
  series} of $\Gamma$ is the formal power series
$$f_\Gamma(X)=\sum_{w\in L(\Gamma)}X^{|w|}.$$

As examples, consider $\Gamma=(\{S\},\{a\},S,\{S\to a,S\to aaS\})$.
Then $L(\Gamma)$ consists precisely of those words $a^n$ for which $n$
is odd. Indeed the only derivation path is $S\to aaS
\Longrightarrow\cdots \Longrightarrow a^{2n}S \Longrightarrow
a^{2n+1}$. This grammar is right-linear and unambiguous. Consider also
$\Gamma'=(\{S\},\{(,)\},S,\{S\to(S),S\to SS,S\to\varepsilon\})$, where
$\varepsilon$ denotes the empty sequence.
$L(\Gamma')=\{\varepsilon,(),()(),(()),\dots\}$ consists precisely of
those words over $\{(,)\}$ that correspond to legal bracket nestings.
This grammar is neither left- nor right-linear.

\subsection{Ergodic regular languages}
Given a context-free grammar $\Gamma = (N,T,S,R)$, its
\emph{dependency di-graph}~\cite{kuich:entropy} $\DD = \DD(\CC)$ is the
oriented graph with vertex set $N$, with an edge from $X$ to $Y$
(written $X \longrightarrow Y$) whenever there is a production $X \to y$ with the
non-terminal $Y$ appearing in $y\in(N\cup T)^*$. We write $X \tos Y$ if in $\DD$
there is an oriented path of length $\ge 0$ from $X$ to $Y$. Recall
that a subset $N'\subseteq N$ is \emph{strongly connected} if for all
$X, Y \in N'$ one has $X \tos Y$ and $Y \tos X$ in $\DD$.

The start symbol $S$ of $\Gamma$ is \emph{isolated} if it does not
occur in the right hand side of any production rule.

%For any language $L \subset T^*$, we define its {\it subword closure\/} 
%as $Sub(L) = \{ v \in T^* : v \;\text{is a subword of}\; w 
%\;\text{for some}\; w \in L \}$.

\begin{defn}
  A linear grammar $\Gamma$ is \emph{ergodic} if (1) either $N$ or
  $N\setminus\{S\}$ is strongly connected, and (2) every terminal word
  $w \in T^*$ occurs as a subword of $v$ in a intermediate derivation
  $U \to v$ with $U$ respectively in $N$ or $N\setminus\{S\}$ and
  $v\in(N \cup T)^*$).
  
  A regular language $L$ is ergodic if it is generated by an ergodic,
  reduced right-linear grammar.
\end{defn}
Note that if $N\setminus\{S\}$ is strongly connected, then clearly $S$
is isolated.  See~\cite{ceccherini-w:ergodicity} for more on
ergodicity of context-free grammars.

\subsection{Growth series and systems of equations associated with grammars}
We shall use the following results:
\begin{thm}[Chomsky and Sch\"utzenberger~\cite{chomsky-s:alg}]\label{thm:ch}
  Let $\Gamma$ be an unambiguous context-free grammar. Then its growth
  series is an algebraic function, i.e.\ there is a polynomial $0\neq
  P(X,Y)\in\Z[X,Y]$ such that $P(X,f_\Gamma(X))=0$.
\end{thm}

\begin{thm}[\cite{salomaa-s:fps,harrison:flt}]\label{thm:reg}
  Let $\Gamma$ a left- or right-linear grammar. Then its growth series
  is rational, i.e.\ there exist polynomials $P,Q\in\Z[X]$ with
  $f_\Gamma(X)=P(X)/Q(X)$.
\end{thm}

\begin{proof}[Sketch of proof] The idea of the proof is to transform a
  context-free grammar into a system of algebraic equations for the
  associated formal power series. Let $\Gamma=(N,T,S,R)$ be any
  reduced (i.e. without superfluous variables), context-free grammar
  without chain and $\epsilon$-rules (i.e.\ productions of the form $U
  \to V$ and $U \to \epsilon$ respectively, with $U,V\in N$). With
  each $U\in N$ we associate the generating function (power series)
  $$f_U(X) = \sum_{w \in L(\Gamma), U \Derives w}X^{|w|}.$$
  In
  addition to $X$, consider extra commuting variables $y_U$ for $U \in
  N$.  Define $\pi(a) = X$ for every $a \in T$ and $\pi(U) = y_U$ for
  every $U \in N$, and for $u=u_1 \cdots u_k \in (T \cup N)^*$, set
  $\pi(u) = \pi(u_1)\cdots\pi(u_k)$, the product of the corresponding
  variables. With $U \in N$ we associate the polynomial
  $${\mathcal P}_U(z; \{y_V\}_{V\in N}) = \sum_{U \to u} \pi(u)$$
  in the variables $X$ and $\{y_U\}$. Then the functions $f_U(X)$
  satisfy the system of equations
  \begin{equation}\label{Chomsky}
   \bigg\{f_U(X) = {\mathcal P}_U\bigl( X; \{f_V(X)\}_{V\in N}\bigr)\bigg\},\quad (U \in N). 
  \end{equation}
  Since each $f_U(X)$ is a power series with non-negative coefficients,
  its radius of convergence $r_U$ is the smallest positive singularity
  of $f_U(X)$, by Pringsheim's theorem. Denoting by $\gamma(L) =
  \limsup_{n \to \infty} \sqrt[n]{|\set{w \in L}{\ell(w) = n}|}$ the
  exponential growth rate of a language $L$, we have
  $$\gamma(L(\Gamma)) = 1/{r_S}.$$
  
  Finally, if $U,V\in N$ and $U\tos V$ in $\DD(\CC)$, then $r_U\le
  r_V$. In particular, if the grammar is ergodic, the radius of
  convergence is the same for all variables: $r_S=r_U=r$ for all $U\in
  N$.
\end{proof}

\subsection{Growth of vertices in $\xlm$}
We identify the set of vertices in $\xlm$ with a language over
$\{e_1,\dots,e_\ell,e_p\}$, by selecting for each vertex $x\in
V(\xlm)$ the label of the unique left-most geodesic from $*$ to $x$.
This language can be described by a right-linear grammar, which we
first present for even $m=2w$.
\begin{itemize}
\item its set of non-terminals is $N = \{X_0, X_1, X_{L,2}, \ldots,
  X_{L,w}, X_{R,2},\ldots,X_{R,w-1}\}$;
\item its set of terminals is $T=\{e_1,\dots,e_\ell\}$;
\item its axiom is $X_0$;
\item its rules are
  \begin{align*}
   X_0 &\to \epsilon | e_1X_1 | \dots | e_\ell X_1\\
   X_1 &\to \epsilon | e_1 X_{R,2} | e_2 X_1 | \dots |
    e_{\ell-2}X_1 | e_{\ell-1}X_{L,2}\\
   X_{L,t} &\to \epsilon | e_1 X_{R,2} | e_2 X_1 | \dots |
    e_{\ell-2}X_1 | e_{\ell-1}X_{L,t+1}\text{ for }2\le t\le w-1\\
   X_{R,t} &\to \epsilon | e_1 X_{R,t+1} | e_2 X_1 | \dots |
    e_{\ell-2}X_1 | e_{\ell-1}X_{L,2}\text{ for }2\le t\le w-2\\
   X_{R,w-1} &\to \epsilon | e_2 X_1 | \dots | e_{\ell-2}X_1 |
    e_{\ell-1}X_{L,2}\\
   X_{L,w} &\to \epsilon | e_1X_{R,2} | e_2X_1 | \dots |
    e_{\ell-3}X_1 | e_{\ell-2}X_{L,2}.
  \end{align*}
\end{itemize}

\begin{prop}
  The language generated by the above grammar is ergodic and it is a
  geodesic normal form for the vertices of $\xlm$.
\end{prop}
\begin{proof}
  The non-terminal $X_0$ expresses the cone of $*$, that is, the whole
  graph; $X_1$ expresses the cone of a vertex of type $1$. The
  non-terminal $X_{L,t}$ expresses the cone of a type-$t$ vertex on the
  left of a $2$-cell, and similarly for $X_{R,t}$.

  The set of rules for a non-terminal $X_t$, $X_{L,t}$ or $X_{R,t}$
  corresponds to the decomposition of the cone of a vertex of extended cone
  type $t, (L,t)$ or $(R,t)$
  into subcones at its successors.
  
  The vertices in the cone at a vertex $x=e_{i_1}\dots e_{i_\ell}$ are
  derived from an intermediate derivation of the form $e_{i_1}\dots
  e_{i_\ell}X_s$, for some
  $s\in\{0,1,(L,2),\dots,(L,w-1),(R,2),\dots,(L,w)\}$. For instance,
  $x$ itself is derived from $xX_s$ by the rule $X_s\to\epsilon$.

  Note that the only difference between left and right is the
  extra rule $X_{L,w-1}\to e_{\ell-1}X_{L,w}$ which has no right
  counterpart. In that way the vertex of type $w$ derived from
  $X_{L,w}$, which has two predecessors, is produced only once from
  its left predecessor.
  
  Ergodicity is obvious.
\end{proof}

We now present the (ergodic) grammar for $\V(\xlm)$ for odd $m=2w+1$,
omitting its analogous proof:
\begin{itemize}
\item its set of non-terminals is $N = \{X_0, X_1, X_{L,2}, \ldots,
  X_{L,m}, X_{R,2},\ldots,X_{R,m-1}\}$;
\item its set of terminals is $T=\{e_1,\dots,e_\ell\}$;
\item its axiom is $X_0$;
\item its rules are
  \begin{align*}
   X_0 &\to \epsilon | e_1X_1 | \dots | e_\ell X_1\\
   X_1 &\to \epsilon | e_1X_{R,2} | e_2 X_1 | \dots |
    e_{\ell-2}X_1 | e_{\ell-1}X_{L,2}\\
   X_{L,t} &\to \epsilon | e_1 X_{R,2} | e_2 X_1 | \dots |
    e_{\ell-2}X_1 | e_{\ell-1}X_{L,t+1}\text{ for }2\le t\le m-1,
    t\neq w\\
   X_{R,t} &\to \epsilon | e_1 X_{R,t+1} | e_2 X_1 | \dots |
    e_{\ell-2}X_1 | e_{\ell-1}X_{L,2}\text{ for }2\le t\le m-2,
    t\neq w\\
   X_{L,w} &\to \epsilon | e_1 X_{R,2} | e_2 X_1 | \dots |
    e_{\ell-3}X_1 | e_{\ell-2}X_{L,t+1}\\
   X_{R,w} &\to \epsilon | e_1 X_{R,t+1} | e_2 X_1 | \dots |
    e_{\ell-3}X_1 | e_{\ell-2}X_{L,2}\\
   X_{R,m-1} &\to \epsilon | e_2 X_1 | \dots | e_{\ell-2}X_1 |
    e_{\ell-1}X_{L,2}\\
   X_{L,m} &\to \epsilon | e_1X_{R,2} | e_2X_1 | \dots |
    e_{\ell-3}X_1 | e_{\ell-2}X_{L,2}.
  \end{align*}
\end{itemize}

Solving the corresponding system~(\ref{Chomsky}) one has:
\begin{cor} The growth series $F_{\ell,m}$ of the vertices in $\xlm$
  with even $m=2w$ is
  \[F_{\ell,m}(X) = \frac{1 + 2X + \dots + 2X^{w-1} + X^w}{ 1 -
   (\ell-2)X - \dots - (\ell-2)X^{w-1} + X^w};\]
  for odd $m=2w+1$, it is
  $$F_{\ell,m}(X) = \frac{1 + 2X + \dots + 2X^{w-1} + 4X^w + 2X^{w+1} +
  \dots + 2X^{m-2} + X^{m-1}}{
  1 - (\ell-2)X - \dots - (\ell-4)X^w - \dots - (\ell-2)X^{m-2} + X^{m-1}}.$$
\end{cor}

\subsection{Growth of geodesics in $\xlm$}\label{sec:ggeod}
We now compute the growth series associated to finite geodesics
starting at $*$; that is,
$$G_{\ell,m}(X) = \sum_{\gamma\in\geod_*}X^{|\gamma|}.$$

The grammar for the set of geodesics can be obtained from the previous
ones just by adding the extra rules
$$X_{R,w-1} \to x_1 X_{L,w}$$
for even $m$, and
$$X_{R,m-1} \to x_1 X_{L,m}$$
for odd $m$. Indeed one can reach a
vertex of type $w$ (resp. of type $m$) both from the right and the
left. While in the grammar for vertices such a vertex is derived only
from the left-side (otherwise we would count it twice!) here the two
geodesic ending at this point are distinct. This grammar is also ergodic.

Solving the associated system~(\ref{Chomsky}) one has:
\begin{thm} The growth series $G_{\ell,m}$ of the geodesics in $\xlm$
with even $m=2w$ is
$$G_{\ell,m}(X) = \frac{1 + 2X + \dots + 2X^{w-1} + X^w + (X^w - X^{w-1})}{
  1 - (\ell-2)X - \dots - (\ell-2)X^{w-1} + X^w + (X^w - X^{w-1})};$$
for odd $m=2w+1$ it is
$$G_{\ell,m}(X) = \frac{\multifrac{1 + 2X + \dots + 2X^{w-1} + 4X^w
   +}{\hspace{5em}{} + 2X^{w+1} + \dots + 2X^{m-2} + X^{m-1} + (X^{m-1} - X^{m-2})}}{
  \multifrac{1 - (\ell-2)X - \dots - (\ell-4)X^w -}{{} - (\ell-2)X^{w+1} -
\dots - (\ell-2)X^{m-2} + X^{m-1} + (X^{m-1} - X^{m-2})}}.$$
\end{thm}

\subsection{Growth of ordered pairs of geodesics in $\xlm$}\label{subs:pairs}
We compute the growth series of ordered pairs of geodesics, both
starting at the base point $*\in V(\gf)$ and having the same endpoint;
that is,
$$H_{\ell,m}(X) =
\sum_{\gamma,\delta\in\geod_*,\gamma_{|\gamma|}=\delta_{|\delta|}}X^{|\gamma|}.$$
We express this set of pairs of geodesics as a language in $T^*\times
T^*$, that is, a direct product of free monoids. This means that, for
instance,
$(e_{i_1},e_{i_2})(e_{i_3},e_{i_4})=(e_{i_1}e_{i_3},e_{i_2}e_{i_4})$,
where $e_{i_j} \in T$.

We start with even $m=2w$:
\begin{itemize}
\item its set of non-terminals is $$N = \{X_0, X_1, X_{L,2}, \ldots,
  X_{L,w}, X_{R,2},\ldots,X_{R,w}, X_{L,\wedge}, X_{R,\wedge}\};$$
\item its set of terminals is $T\times T$, with
  $T=\{e_1,\dots,e_\ell\}$;
\item its axiom is $X_0$;
\item its rules, with $2\le t\le w-1$, are
\begin{align*}
  X_0\to & \epsilon | (e_1,e_1)X_1 | \dots | (e_\ell,e_\ell)X_1 |\\
      & (e_1e_{\ell-1}^{w-2},e_{2}e_1^{w-2})X_{L,\wedge} | \dots
      | (e_{\ell-2}e_{\ell-1}^{w-2},e_{\ell-1}e_1^{w-2})X_{L,\wedge} |\\
      & (e_{2}e_1^{w-2},e_1e_{\ell-1}^{w-2})X_{R,\wedge} | \dots
      | (e_{\ell-1}e_1^{w-2},e_{\ell-2}e_{\ell-1}^{w-2})X_{R,\wedge}\\
  X_1\to & \epsilon | (e_1,e_1)X_{R,2} | (e_2,e_2)X_1 | \dots
      | (e_{\ell-2},e_{\ell-2})X_1 | (e_{\ell-1},e_{\ell-1})X_{L,2} |\\
      & (e_1e_{\ell-1}^{w-2},e_{2}e_1^{w-2})X_{L,\wedge} | \dots
      | (e_{\ell-2}e_{\ell-1}^{w-2},e_{\ell-1}e_1^{w-2})X_{L,\wedge} |\\
      & \to(e_{2}e_1^{w-2},e_1e_{\ell-1}^{w-2})X_{R,\wedge} | \dots
      | (e_{\ell-1}e_{1}^{w-2},e_{\ell-2}e_{\ell-1}^{w-2})X_{R,\wedge}\\
\end{align*}
\begin{align*}
  X_{L,t}\to & \epsilon | (e_1,e_1)X_{R,2} | (e_2,e_2)X_1 | \dots
      | (e_{\ell-2},e_{\ell-2})X_1 | (e_{\ell-1},e_{\ell-1})X_{L,t+1} |\\
      & (e_1e_{\ell-1}^{w-2},e_{2}e_1^{w-2})X_{L,\wedge} | \dots
      | (e_{\ell-2}e_{\ell-1}^{w-2},e_{\ell-1}e_1^{w-2})X_{L,\wedge} |\\
      & (e_{2}e_1^{w-2},e_1e_{\ell-1}^{w-2})X_{R,\wedge} | \dots
      | (e_{\ell-1}e_1^{w-2},e_{\ell-2}e_{\ell-1}^{w-2})X_{R,\wedge}\\
  X_{R,t}\to & \epsilon | (e_1,e_1)X_{R,t+1} | (e_2,e_2)X_1 | \dots
      | (e_{\ell-2},e_{\ell-2})X_1 | (e_{\ell-1},e_{\ell-1})X_{L,2} |\\
      & \to(e_1e_{\ell-1}^{w-2},e_{2}e_1^{w-2})X_{L,\wedge} | \dots
      | (e_{\ell-2}e_{\ell-1}^{w-2},e_{\ell-1}e_1^{w-2})X_{L,\wedge} |\\
      & (e_{2}e_1^{w-2},e_1e_{\ell-1}^{w-2})X_{R,\wedge} | \dots
      | (e_{\ell-1}e_1^{w-2},e_{\ell-2}e_{\ell-1}^{w-2})X_{R,\wedge}\\
  X_{R,w}\to & \epsilon | (e_1,e_1)X_{R,2} | (e_2,e_2)X_1 | \dots
      | (e_{\ell-3},e_{\ell-3})X_1 | (e_{\ell-2},e_{\ell-2})X_{L,2} |\\
      & \to(e_1e_{\ell-1}^{w-2},e_{2}e_1^{w-2})X_{L,\wedge} | \dots
      | (e_{\ell-3}e_{\ell-1}^{w-2},e_{\ell-2}e_1^{w-2})X_{L,\wedge} |\\
      & \to(e_{2}e_1^{w-2},e_1e_{\ell-1}^{w-2})X_{R,\wedge} | \dots
      | (e_{\ell-2}e_1^{w-2},e_{\ell-3}e_{\ell-1}^{w-2})X_{R,\wedge}\\
  X_{L,w}\to & X_{R,w}\\
  X_{L,\wedge}\to & (e_{\ell-1},e_1)X_{L,w}
      | (e_{\ell-1}^{w-1},e_2e_1^{w-1})X_{L,\wedge}
      | (e_{\ell-2}e_{\ell-1}^{w-2},e_1^{w-1})X_{L,\wedge}\\
  X_{R,\wedge}\to & (e_1,e_{\ell-1})X_{R,w}
      | (e_2e_1^{w-1},e_{\ell-1}^{w-1})X_{R,\wedge}
      | (e_1^{w-1},e_{\ell-2}e_{\ell-1}^{w-2})X_{R,\wedge}\\
\end{align*}
\end{itemize}

\begin{prop}\label{prop:grampair}
  The language generated by the above grammar is ergodic and it is a 
  geodesic normal form for the ordered pairs of
  geodesics with common extremities.
\end{prop}
\begin{proof}
  Given a geodesic $\gamma=(\gamma_0,\dots,\gamma_n)$ the truncation
  of $\gamma$ is $\gamma' =(\gamma_0,\dots,\gamma_{n-1})$; more
  generally, for $i\le n$, define its $i$-step truncation
  $\gamma^{(i)}= (\gamma^{(i-1)})' = (\gamma_0,\dots,\gamma_{n-i})$,
  the geodesic obtained by deleting its last $i$ edges.

  We first explain the meaning of the non-terminals. If
  $(\gamma,\delta)$ ends in a vertex of extended type $t$, then
  $X_0\Derives(\gamma,\delta)X_t$ is an intermediate derivation.

  Therefore $X_{L,t}$ expresses the possible continuations of a pair
  of geodesics that pass through a common point of type $t$ on the
  left of a $2$-cell --- and analogously for $X_{R,t}$.

  If $(\gamma,\delta)$ is a pair of geodesics of length $n$ with
  $\gamma_{n-1}\neq\delta_{n-1}$, then (by Lemma~\ref{lem:hyper})
  there is a unique $2$-cell $F$ below $\gamma_n$ and between $\gamma$
  and $\delta$. If $\gamma$ is on the left and $\delta$ is on the
  right of $F$, then $X_0\Derives(\gamma',\delta')X_{L,\wedge}$ is an
  intermediate derivation. If $\gamma$ is on the right and $\delta$ is
  on the left of $F$, then $X_0\Derives(\gamma',\delta')X_{R,\wedge}$
  is an intermediate derivation.

  Therefore $X_{L,\wedge}$ and $X_{R,\wedge}$ express the possible
  continuations of the $1$-step truncation of a pair of geodesics
  surrounding a $2$-cell.

  Let $(\gamma,\delta)$ be a pair of geodesics of length $n$. We show
  by induction on $n$ that $(\gamma,\delta)$ is uniquely derived by
  the above grammar. Write $t$ the extended type of $\gamma_n$. We
  distinguish several cases:
  \begin{description}
  \item[if $\mathbf{\gamma_{n-1}=\delta_{n-1}}$] then there is a unique
   derivation
   \[X_0\Derives(\gamma',\delta')X_{t'} \Longrightarrow(\gamma,\delta)X_t
    \Longrightarrow(\gamma,\delta),\]
   where $t'$ is the extended type of $\gamma_{n-1}$;
  \item[if $\mathbf{\gamma_{n-1}\neq\delta_{n-1}}$] let $i$ be minimal
   such that $\gamma_{n-1-i(w-1)}=\delta_{n-1-i(w-1)}$. Then there is
   a unique derivation
   \begin{multline*}
   X_0\Derives(\gamma^{(1+i(w-1))},\delta^{(1+i(w-1))})X_{t'}
    \Longrightarrow(\gamma^{(1+(i-1)(w-1))},\delta^{(1+(i-1)(w-1))})X_{S,\wedge}\\
    \Longrightarrow\dots \Longrightarrow(\gamma',\delta')X_{S,\wedge}
    \Longrightarrow(\gamma,\delta)X_t \Longrightarrow(\gamma,\delta),
   \end{multline*}
   where $t'$ is the extended type of $\gamma_{n-1-i(w-1)}$ and $S\in\{L,R\}$
   determines which of $\gamma$ and $\delta$ is on the left between
  $\gamma_{n-1-i(w-1)}$ and $\gamma_n$.
  \end{description}
\end{proof}

We now present the (ergodic) grammar for odd $m=2w+1$; we omit the
proof --- completely analogous to that of
Proposition~\ref{prop:grampair} --- that it describes pairs of
geodesics, since
\begin{itemize}
\item its set of non-terminals is $$N = \{X_0, X_1, X_{L,2}, \ldots,
  X_{L,m}, X_{R,2},\ldots,X_{R,m}, X_{L,\wedge}, X_{R,\wedge}\};$$
\item its set of terminals is $T=\{e_1,\dots,e_\ell\}$;
\item its axiom is $X_0$;
\item its rules are
\begin{align*}
 X_0 \to & \epsilon | (e_1 ,e_1) X_1 | \dots |(e_\ell,e_\ell) X_1|\\
      & (e_1 e_{\ell-1}^{m-2}, e_{2} e_1^{m-2})X_{L,\wedge} | \dots
      | (e_{\ell-2} e_{\ell-1}^{m-2},e_{\ell-1}e_1^{m-2})X_{L,\wedge}|\\
      & (e_{2} e_1^{m-2}, e_1e_{\ell-1}^{m-2}) X_{R,\wedge}| \dots
      | (e_{\ell-1} e_1^{m-2}, e_{\ell-1}e_{\ell-1}^{m-2})X_{R,\wedge}\\
 X_1 \to & \epsilon | (e_1,e_1)X_{R,2} | (e_2,e_2) X_1 | \dots
      | (e_{\ell-2},e_{\ell-2})X_1 | (e_{\ell-1},e_{\ell-1})X_{L,2} |\\
      & (e_1 e_{\ell-1}^{m-2},e_{2} e_1^{m-2})X_{L,\wedge} | \dots
      | (e_{\ell-2}e_{\ell-1}^{m-2},e_{\ell-1}e_1^{m-2})X_{L,\wedge} |\\
      & \to(e_{2}e_1^{m-2},e_1e_{\ell-1}^{m-2})X_{R,\wedge}| \dots
      | (e_{\ell-1} e_{1}^{m-2},e_{\ell-2}e_{\ell-1}^{m-2})X_{R,\wedge}\\
X_{L,t}\to & \epsilon | (e_1,e_1)X_{R,2} | (e_2,e_2)X_1 | \dots
      | (e_{\ell-2},e_{\ell-2})X_1 | (e_{\ell-1},e_{\ell-1})X_{L,t+1} |\\
      & (e_1e_{\ell-1}^{m-2},e_{2}e_1^{m-2})X_{L,\wedge} | \dots
      | (e_{\ell-2}e_{\ell-1}^{m-2},e_{\ell-1}e_1^{m-2})X_{L,\wedge} |\\
      & (e_{2}e_1^{m-2},e_1e_{\ell-1}^{m-2})X_{R,\wedge} | \dots
      | (e_{\ell-1}e_1^{m-2},e_{\ell-2}e_{\ell-1}^{m-2})X_{R,\wedge}\\
X_{R,t}\to & \epsilon | (e_1,e_1)X_{R,t+1} | (e_2,e_2)X_1 | \dots
      | (e_{\ell-2},e_{\ell-2})X_1 | (e_{\ell-1},e_{\ell-1})X_{L,2} |\\
      & \to(e_1e_{\ell-1}^{m-2},e_{2}e_1^{m-2})X_{L,\wedge} | \dots
      | (e_{\ell-2}e_{\ell-1}^{m-2},e_{\ell-1}e_1^{m-2})X_{L,\wedge} |\\
      & (e_{2}e_1^{m-2},e_1e_{\ell-1}^{m-2})X_{R,\wedge} | \dots
      | (e_{\ell-1}e_1^{m-2},e_{\ell-2}e_{\ell-1}^{m-2})X_{R,\wedge}\\
X_{R,w}\to & \epsilon | (e_1,e_1)X_{R,w+1} | (e_2,e_2)X_1 |\dots
      | (e_{\ell-3},e_{\ell-3}) X_1| (e_{\ell -2},e_{\ell-2}) X_{L,2} |\\
      & \to (e_1e_{\ell-1}^{m-2},e_{2}e_1^{m-2})X_{L,\wedge} | \dots
      | (e_{\ell-3}e_{\ell-1}^{m-2},e_{\ell-2} e_1^{m-2}) X_{L,\wedge} |\\
      & (e_{2}e_1^{m-2},e_1e_{\ell-1}^{m-2})X_{R,\wedge} | \dots 
      | (e_{\ell-2}e_1^{m-2},e_{\ell-3}e_{\ell-1}^{m-2})X_{R,\wedge}\\
X_{L,w}\to & \epsilon | (e_1,e_1)X_{R,2} |(e_2,e_2) X_1 |\dots
      | (e_{\ell-3},e_{\ell-3})X_1| (e_{\ell -2},e_{\ell-2})X_{L,w+1} |\\
      & (e_1e_{\ell-1}^{w-2},e_{2}e_1^{w-2})X_{L,\wedge} |\dots 
      | (e_{\ell-3}e_{\ell-1}^{w-2},e_{\ell-2}e_1^{w-2})X_{L,\wedge} |\\
      & (e_{2}e_1^{w-2},e_1e_{\ell-1}^{w-2})X_{R,\wedge}| \dots
      | (e_{\ell-2} e_1^{w-2},e_{\ell-3}e_{\ell-1}^{w-2})X_{R,\wedge}\\
X_{R,m} \to & \epsilon | (e_1,e_1)X_{R,2} | (e_2 ,e_2) X_1 |\dots
      | (e_{\ell-3},e_{\ell-3})X_1| (e_{\ell -2},e_{\ell-2})X_{L,2} |\\
      & (e_1e_{\ell-1}^{m-2},e_{2}e_1^{m-2})X_{L,\wedge} |\dots 
      | (e_{\ell-3}e_{\ell-1}^{m-2},e_{\ell-2}e_1^{m-2})X_{L,\wedge} |\\
      & (e_{2}e_1^{m-2},e_1e_{\ell-1}^{m-2})X_{R,\wedge} | \dots 
      | (e_{\ell-2}e_1^{m-2},e_{\ell-3}e_{\ell-1}^{m-2})X_{R,\wedge}\\
X_{L,m} \to & X_{R,m}\\
X_{L,\wedge} \to & (e_{\ell-1},e_1)X_{L,m} |
              (e_{\ell-1}^{m-1},e_2e_1^{m-1})X_{L,\wedge} |
              (e_{\ell-2}e_{\ell-1}^{m-2},e_1^{m-1})X_{L,\wedge}\\
X_{R,\wedge} \to & (e_1,e_{\ell-1})X_{R,m} |
              (e_2e_1^{m-1},e_{\ell-1}^{m-1})X_{R,\wedge} |
              (e_1^{m-1},e_{\ell-2}e_{\ell-1}^{m-2})X_{R,\wedge}
\end{align*}
\end{itemize} 

\begin{thm} The growth series $H_{\ell,m}$ of the pair of geodesics in
    $\xlm$ with even $m=2w$ is
$$H_{\ell,m}(X) = 1 + \ell\frac{X + X^2 + \dots + X^{w-1}}{
  \multifrac{1 - (\ell-2)X - \dots - (\ell-2)X^{w-1} + X^w +}{{} +
X^{w-1}(X-1)(2X^{w-1} - 3)}};$$
for odd $m=2w+1$, it is
$$H_{\ell,m}(X) = 1 + \ell\frac{X + X^2 + \dots + X^{m-2}}{
  \multifrac{1 - (\ell-2)X - \dots - (\ell-4)X^{w-1} - \dots -
  (\ell-2)X^{m-2} +}{{} + X^{m-1} - X^{m-1}(1-X)(3-4X^w+2X^{m-2})}}.$$
\end{thm}

\subsection{Holly trees}\label{subs:holly}
Regular grammars were used in the previous subsection to describe
fellow-traveling pairs of geodesics, giving in the next section a
lower estimation on the cogrowth (and hence the spectral radius) of
$\xlm$, using loops of the form $\gamma\rho\delta^{-1}$, where
$(\gamma,\delta)$ is a pair of geodesics with common extremities $*$
and $x$, and $\rho$ is the loop around a $2$-cell at $x$, not tangent
to $\gamma$ nor $\delta$.

We give here a context-free grammar producing a larger class of
 proper loops, which we call \emph{holly trees}. To define them, call
$\mathcal H_x$ the set of holly trees at the vertex $x\in\V(\xlm)$; then
$\mathcal H_x$ is a set of loops based at $x$, and lying entirely in
the cone of $x$. Then
\begin{enumerate}
\item if $\gamma,\delta\in\mathcal H_x$ and $\gamma\delta$ is a
  proper path, then $\gamma\delta\in\mathcal H_x$;
\item if $y$ is a successor of $x$ and $\gamma\in\mathcal H_y$, then
  $v\gamma v^{-1}\in\mathcal H_x$, where $v$ is the label of the edge
  from $x$ to $y$;
\item if $P$ is a $2$-cell in the cone of $x$ and touching $x$, whose
  perimeter starting at $x$ is labelled $\gamma$, then
  $\gamma\in\mathcal H_x$;
\item no other loop belongs to $\mathcal H_x$.
\end{enumerate}

We note that the growth of $\mathcal H_x$ depends only on the cone
type of $x$, since all holly trees at $x$ lie within the cone of $x$. We
may thus consider $L_t(X)$, the growth series of holly trees at any
fixed vertex of type $t$. That this function is algebraic follows from
the fact that holly trees can be described by the unambiguous
context-free grammar given below (and compare with Theorem~\ref{thm:ch}).

Note also that all holly-trees are non-trivial proper paths. Write
$\epsilon$ for the empty path.

For simplicity, assume $m=2w$ is even. Consider the non-terminal $L_{t,e}$
for all cone types $t$, and all labels $e$ of the successors of a
vertex of type $t$. We then have $t\in\{0,\dots,w\}$ and
$e\in\{e_1,\dots,e_{\ell(t)}\}$. The
variable $L_{t,e_i}$ expresses the set of all holly trees in the cone (at
a vertex of) type $t$,
whose first edge is labelled $e_i$. For commodity, define also
variables $L_t$, counting all holly trees in a cone of type $t$, and
$L_{t,\hat e}$, expressing all holly trees in a cone of type $t$ whose
first edge is not labelled $e$.

The terminal alphabet is $\{e_1^{\pm1},\dots,e_\ell^{\pm1}\}$,
describing paths in $\xlm$. The axiom is $L_0$. The derivations are
\begin{align*}
  L_t &\to L_{t,e_1} | \cdots | L_{t,e_{\ell(t)}}\text{ for all }t\in\{0,\dots,w\}\\
  L_{t,\hat e} &\to\text{the same as $L_t$, but excluding }L_{t,e}\\
  L_{t,e} &\to eL_ue^{-1}|eL_ue^{-1}L_{t,\hat e}| e_1\dots e_m |
   e_1\dots e_mL_{t,\widehat{e_m^{-1}}},
\end{align*}
where $e$ leads to a vertex of type $u$ and $(e=e_1,\dots,e_m)$
describes all possible labelings starting with $e$ around $2$-cells in
a cone of type $t$.

The computations of growth series are more complicated, and the
results, for $m=\ell=8$, appear in the next section.

\begin{remark}
  We end this section by remarking that also Cayley graphs of
  hyperbolic groups, or, more generally strongly-transitive hyperbolic
  graphs have a finite number of cones types. As a consequence,
  \emph{mutatis mutandis}, all our above computations (enumeration of
  vertices; of geodesics; of ordered pairs of geodesic or, more
  generally, of ordered $N$-tuples of geodesics; of holly trees) can
  be performed for these groups as well (always leading to rational
  growth series, except for the holly trees for which the growth
  series will again be algebraic).
\end{remark}

\section{Estimates for Simple Random Walks}\label{sec:lower}
\subsection{Upper estimates} Upper estimates for the spectral radius
of the Markov operator associated with a simple random walk on the
fundamental group of an orientable surface have been obtained
in~\cite{b-c-c-h:growth}. In that paper various methods were described,
yielding as best estimate for $\ell=m=4g$
$$\mu_g \le \frac{\sqrt{4g-2}}{2g}+\frac1{4g},$$
and in particular $$\mu_2 \le 0.7374.$$

These estimates (again for $\ell=m=4g$) have been improved by
\.Zuk~\cite{zuk:norm} and by Nagnibeda~\cite{nagnibeda:upperbd}, who
obtained the best known estimate from above:
$$\mu_2 \le 0.6629.$$
Both methods easily extend to all graphs $\xlm$.

\subsection{Lower estimates} In this section, using the Grigorchuk
formula~(\ref{eq:grigorchuk}) and the combinatorial results of
Section~\ref{sec:comput}, we obtain tighter and tighter lower
estimates for $\mu_{\ell,m}$, the spectral radius of the simple random
walk on $\xlm$, and present numerical results for $\mu_{8,8} = \mu_2$.

We start by re-obtaining Kesten's estimate. Consider the loops
consisting of a proper path $\gamma =(*, \ldots, x)$ starting at $*$,
followed by the boundary of a cell containing $x$ and inside its cone,
followed by $\gamma^{-1}$. There are $\ell(\ell-1)^{n-1}$ paths
$\gamma$ of length $n$, and $2(\ell-2)$ choices for the cell's
boundary; thus there are at least
$$\beta_n = 2(\ell-2)\ell(\ell-1)^{n-1}$$
paths of length $2n+m$, whence we re-obtain Grigorchuk's estimate
$$\alpha\ge\limsup_{n\to\infty}\sqrt[2n+m]{\beta_n} = \sqrt{\ell-1},$$
and Kesten's estimate
$$\mu_{\ell,m}\ge\frac{2\sqrt{\ell-1}}{\ell};$$
in particular
$$\mu_2 \ge 0.66143.$$

Taking into account the boundaries of two cells neighbouring $*$
as other constituents of the loops we count, it is possible to
obtain slightly tighter results. As the gain is negligible, we
shall not describe the counting in detail, but refer to
the paper by Kesten~\cite[Theorem~4.15]{kesten:rwalks}.

The first non-trivial estimate is obtained by considering the growth
of \emph{unimodular loops}. A \emph{unimodular loop} of weight $n$ is
a loop $\gamma = (\gamma_0, \gamma_1, \ldots, \gamma_n, \gamma_{n+1},
\ldots, \gamma_{2n})$ of length $2n$ such that $|\gamma_{n}| = n$ and
such that $\gamma_{n-1} \neq \gamma_{n+1}$. In other words $\gamma$
can be viewed as a couple of distinct geodesics of length $n$, say
$\delta=(\delta_0, \ldots, \delta_n)$ and $\beta=(\beta_0, \ldots,
\beta_n)$ with common starting and end-points, namely $\delta_0 =
\beta_0$ and $\delta_n = \beta_n$ but with no common last edge, namely
$\delta_{n-1} \neq \beta_{n-1}$. We then regard $\gamma$ as
$(\delta_0, \ldots, \delta_n, \beta_{n-1}, \ldots, \beta_0)$.

The growth series $U$ for unimodular loops is closely related to the
function $H_{\ell,m}$ computed in subsection~\ref{subs:pairs}: both
functions have the same denominator, and hence the same radius of
convergence. This is because $U(X) < H_{\ell,m}(X)$ and $H_{\ell,m}(X)
< X^mU(X)$ coefficient-wise: any unimodular loop comes from a pair of
geodesics, and any pair of geodesics can be extended (if necessary) in
at most $m$ steps in a unimodular loop. For $\ell=m=8$, for instance,
the radius of convergence of $U$ is
$$\rho_{U} \approx 1/7.0248,$$
so that $\alpha_2\ge \sqrt{\frac1{\rho_{f_U}}} \approx 2.65$, and thus,
by Remark~\ref{rem:alphamu} and Grigorchuk's formula~(\ref{eq:grigorchuk}),
$$\mu_2\ge 0.66144.$$

This last estimate is weaker than that obtained by
Paschke~\cite[Theorem~3.2]{paschke:norm}, whose result applied to
$\x_{8,8}$ gives
$$\mu_2\ge 0.6616.$$

A sharper estimate is given by the growth of holly trees as computed
in subsection~\ref{subs:holly}. Again in the case $\ell=m=8$, the
algebraic growth function was computed using the computer algebra
program \textsc{Maple}, as a solution $L(X)$ of $P(X,L(X))\equiv 0$.
The algebraic function $L$ has a singularity at all vanishing points
$\rho$ of the discriminant of $P(X,Y)$, so the radius of convergence
of $L$ is at most the absolute value of a minimal such $\rho$. The
discriminant turns out to be a degree-$179$ polynomial vanishing at
$$\rho \approx 0.12887$$
from which we deduce the asymptotic growth of loops is at least
$$\alpha \ge \frac1{\sqrt\rho} \approx 2.7856,$$
so by
Remark~\ref{rem:alphamu} and Grigorchuk's
formula~(\ref{eq:grigorchuk}) we get the following
\begin{thm} The spectral radius $\mu_2$ associated with a simple
  random walk on the fundamental group $J_2$ of a surface of genus $2$
  (with respect to the canonical set of generators) is bounded below
  by
  $$\mu_2 \ge 0.6623.$$
\end{thm}

\begin{remark} 
  The isoperimetric constants $\iota$ for the graphs $\xlm$ have been
  recently computed independently by Yusuke Higuchi and Tomoyuki
  Shirai~\cite{higuchi-s:dfregular} and by Olle H\"aggstr\"om, Johan
  Jonasson and Russell Lyons~\cite{haggstrom-j-l:isoperimetric}. They
  obtained the values
  \begin{equation}
    \iota(\xlm) = (\ell-2)\sqrt{1 - \frac{4}{(\ell - 2)(m - 2)}}.
  \end{equation}
  
  In particular, for $\ell = m = 8$, one has
  $$ \iota(\cay(J_2,S_2)) = \iota(\x_{8,8}) = 4 \sqrt{2}$$
  which, together with (\ref{muiota}) gives
  $$0.38 \approx 1 - \frac{7}{16}\sqrt{2} \leq \mu \leq \frac{1}{\sqrt{2}} \approx 0.7071.$$
  None of these bounds improve our previous estimates.
\end{remark}  

\begin{remark}\label{rem:cactus}
  On the other hand, in~\cite{bartholdi:cactus}, the first author
  introduced a class of loops, called \emph{cactus trees}, which
  contains the class of holly trees and obtained a slightly better
  estimate: $\mu_2 \ge 0.6624$. We do not go into the details here since
  the idea remains the same --- construct an unambiguous context-free
  grammar producing loops --- but the calculations are much more
  intricate.
\end{remark}

\section*{Acknowledgments}
We wish to express our gratitude to Rostislav I.\ Grigorchuk, Pierre
de la Harpe and Wolfgang Woess for many valuable comments and
suggestions.  We also benefitted of stimulating conversions with
Fabrice Liardet and Tatiana Nagnibeda.

The second author wishes to thank the ``Section de Math\'ematiques de
l'Universit\'e de Gen\`eve'' for its kind hospitality during several
stages of the present work; he also acknowledges support from the
C.N.R.\ and from the Swiss National Science Foundation.

This work was completed at the Erwin Schr\"odinger Institute of Vienna
during the conference on ``Random walks and geometry'': we thank the
ESI for its financial support, and the organizers, Vadim Kaimanovich,
Klaus Schmidt and Wolfgang Woess, for their invitation and kind
hospitality.

%\bibliographystyle{amsalpha}
%\bibliography{mrabbrev,people,bartholdi,grigorchuk,math}
\def\nop#1{}\font\cyr=wncyr8\def\cprime{$'$}
\providecommand{\bysame}{\leavevmode\hbox to3em{\hrulefill}\thinspace}

\end{document}